\documentclass{article}
\usepackage{amsmath,amssymb,theorem}

\newcommand{\eop}{\hfill$\square$}
\newcommand\mfk\mathfrak
\newcommand\mcl\mathcal
\newcommand\mbb\mathbb
\newcommand\mtl\mathit
\newcommand\mbf\mathbf
\newcommand\onm\operatorname
\newcommand{\Id}{{\mathchoice
 {\mbox{Id}}{\mbox{Id}}{\mbox{\scriptsize Id}}{\mbox{\scriptsize Id}}}}
\newcommand{\bset}[2]{\bigl\{\,#1:#2\,\bigr\}}
\newcommand{\dbyd}[2]{{\displaystyle\left.\frac d{d#1}\right|_{#1=#2}}}
\newcommand{\semi}{\mathop{_\theta\!\times}}
\newtheorem{lem}{Lemma}
\theorembodyfont{\rmfamily}
\newtheorem{proof}{Proof}

\begin{document}
\title{The Landau--Lifshitz Equation by Semidirect Product Reduction}
\author{George W. Patrick$^*$\\
Department of Mathematics and Statistics\\
University of Saskatchewan\\Saskatoon, Saskatchewan, 
S7N~5E6\\Canada
}
\date{March  2000}
\maketitle
{\renewcommand{\thefootnote}{}
\footnotetext{$^*$Supported in part by the Natural Sciences and Engineering 
Research Council, Canada.}}
\begin{abstract}
The Landau--Lifshitz equation is derived as the reduction of a geodesic
flow on the group of maps into the rotation group. Passing the symmetries
of spatial isotropy to the reduced space is an example of
semidirect product reduction by stages.
\end{abstract}
The Landau--Lifshitz equation for the evolution of a field of unit
vectors $n(x,t)$, $x\in\mbb R^p$, $p=1,2,\cdots$, nondimensionalized, is
\begin{equation*}
\frac{\partial n}{\partial t}=-n\times\frac{\delta E}{\delta n},
\end{equation*}
where $E$ is some energy functional, typically 
\begin{equation*}
E=\frac12\int|\nabla n|^2+\frac a2\int\bigl(n\cdot n-(n\cdot\mbf k)^2\bigr),
\end{equation*}
and $a$ is some coupling constant. This functional encodes the
energy due to orientation variance plus, when $a\ne0$, a material
anisotropy. With these notations $E$ is finite if $n(x)\rightarrow-\mbf k$
sufficiently fast as $x\rightarrow\infty$.
In~\cite{CooperNR-1999.1} and~\cite{PapanicolaouNTomarasTN-1991.1} one
finds, after adjusting dimensions, the conserved quantities
\begin{equation*}
P=2\pi\int x\times\Omega,\qquad\Omega_\alpha=\frac1{8\pi}
 \epsilon_{\alpha\beta\gamma}n\cdot(\nabla_\beta n\times
 \nabla_\gamma n),
\end{equation*}
which generate translations and the conserved
quantity
\begin{equation*}
N=\int(1+n\cdot\mbf k),
\end{equation*}
which generates rotations (uniform in space) of $n$ about $\mbf k$.

In this Letter, I explicitly view the Landau--Lifshitz equation
as a right-hand reduction of a geodesic flow on the group of smooth
mappings
\begin{equation*}
G=\bigl\{A:\mbb R^p\rightarrow \mtl{SO}(3)\bigm|A(\infty)=\Id\bigr\}
\end{equation*}
with pointwise multiplication. The calculations are formal but at the
very least require certain asymptotic behaviors the elements of
$G$. For my purpose it suffices to assume that they are the identity
outside some compact set.

One reason to view the Landau--Lifshitz equation as a reduction from a
cotangent bundle is that, at the level of $T^*G$, the momenta
associated to the group $SE(p)$ and to right translation by $G$ itself
are easily calculated by well-known formulas. Regardless of whether
the energy functional is invariant, one may try to pass these momenta
and their symmetries to the reduced space on which, as it turns out,
the Landau--Lifshitz equation lives. I show that the conserved
quantities $P$ are obtained in this way. The validity of this rests on
whether the subsequent spaces obtained by reductions by those reduced
symmetries are isomorphic to the reductions of $T^*G$ by the whole
symmetry. This is formally correct by the reduction by stages theory
of~\cite{MarsdenJEMisiolekGPerlmutterMRatiuTS-1999}. The earlier
reduction by stages theory found
in~\cite{MarsdenJEMisiolekGPerlmutterM-1998.1} is inadequate for the
Landau--Lifshitz equation since the group $G$ is non-Abelian.

When $p=2$ the field $n$ may be regarded as a map from the $2$-sphere
$S^2$ to itself which maps $\infty$ to $-\mbf k$. Thus, for
$p=2$, the Landau--Lifshitz phase space is decomposed by the
degree of $n$, which is the integer
\begin{equation*}
\onm{deg}n=\frac 1{4\pi}\int n\cdot\frac{\partial n}{\partial x}\times
\frac{\partial n}{\partial y}
\end{equation*}
This normalization is such that the degree is 1 on the identity map
pulled back to the plane using the stereographic projection which
sends $-\mbf k$ to $\infty$. In~\cite{PapanicolaouNTomarasTN-1991.1}
it is observed that the $x$ and $y$ momenta do not commute in
any degree nonzero sector. From a
geometric mechanics point of view, this is an apparent contradiction
since these momenta \emph{do commute} on the unreduced phase
space. Obtaining the degree nonzero sector from the reduction
exercises nontrivial aspects of the reduction by stages
theory and that explains fully the apparent contradiction.

When elements of $G$ are thought of as configurations, I will denote
them as $\psi,\phi\cdots\in G$, while when $G$ plays its role as a
group of symmetries, its elements will be $A,B\cdots\in G$. Formally,
the Lie algebra $\mfk g$ of $G$ is the smooth $\mfk{so}(3)\cong\mbb
R^3$ valued mappings on $\mbb R^p$, vanishing at $\infty$, with
pointwise Lie bracket. The dual $\mfk g^*$ is the same space except
its elements do not necessarily vanish at $\infty$, and the pairing
between $\mfk g$ and $\mfk g^*$ is
\begin{equation*}
\langle\mu,\xi\rangle=\int\mu(x)\cdot\xi(x).
\end{equation*}

I represent the elements of $T^*G$ as pairs $(\psi,\mu)\in G\times\mfk g^*$
using \emph{right} translation, so
the action of right translation of $G$ on $T^*G$ is
\begin{equation*}
A\cdot(\psi,\mu)=\bigl(\psi(x)A(x)^{-1},\mu(x)\bigr).
\end{equation*}
The momentum map of this action is
\begin{equation*}
J^G(\psi,\mu)=-\psi(x)^{-1}\mu(x),
\end{equation*}
because~(see \cite{MarsdenJERatiuTS-1994.1})
\begin{equation*}
\begin{split}
J^G_\xi(\psi,\mu)&=(\psi,\mu)\cdot\dbyd t0\psi\onm{exp}(-\xi t)\\
&=-\int\mu(x)\cdot \psi(x)\xi(x)=\int\bigl(-\psi(x)^{-1}\mu(x)\bigr)
  \cdot\xi(x).
\end{split}\end{equation*}

Fixing the particular momentum value $\mu\in\mfk g^*$ defined by
$\mu(x)=\mbf k$, one has
\begin{equation*}
(J^G)^{-1}(\mu)=\bigl\{(\psi,-\psi\mbf k)\bigr\},\quad 
G_\mu=\bset{A\in G}{A(x)\mbf k=\mbf k}
\end{equation*}
and the quotient $(J^G)^{-1}(\mu)/G_\mu$ may be realized as 
\begin{equation*}
(J^G)^{-1}(\mu)\rightarrow\bset{n\in\mfk g}{|n|=1}\quad\mbox{by}\quad
n=-\psi\mbf k.
\end{equation*}
It is standard that the evolution equations on
$(J^G)^{-1}(\mu)/G_\mu$ (i.e. the evolution equations for $n$)
are the Lie--Poisson (`+' because of the right translation) equations
\begin{equation*}
\frac{\partial n}{\partial t}=+\onm{coad}_{\delta E/\delta n}n
\end{equation*}
which is exactly the Landau--Lifshitz equation. Here, by definition
$\onm{coad}_\xi\mu=-\onm{ad}^*_\xi\mu$, for $\xi\in\mfk g$ and
$\mu\in\mfk g^*$.

The question of whether all fields $n$ are
attained through the reduction is the question of whether any
$n:\mbb R^p\rightarrow\mfk{so}(3)$ lifts through $\psi\mapsto-\psi\mbf
k$ to $\mtl{SO}(3)$ with $\psi(\infty)=\Id$. This is true for any $n$
if $p\ne 2$ and if $\onm{deg}n=0$ in the case $p=2$. Indeed, suppose
$R>0$ is such that $n=-\mbf k$ for $|x|\ge R$. Choosing any connection
on the principle bundle $\mtl{SO}(3)\rightarrow S^2$ by $A\mapsto
-A\mbf k$ and lifting $n$ along the radial lines emanating from the
origin gives a smooth $\tilde\psi(x)$ such that $\tilde \psi(0)=\Id$,
$\tilde\psi$ is constant along radial lines for $|x|\ge R$, and
$\tilde\psi(x)(-\mbf k)=-\mbf k$ for $|x|\ge R$. Thus $\tilde\psi$
restricted to the sphere $S_R^{p-1}(0)$ of radius $R$ centered at $0$
is a map into the isotropy group $\bigl(\mtl{SO}(3)\bigr)_{\mbf
k}\cong S^1$. In the case $p\ne 2$ the homotopy group $\pi_{p-1}(S^1)$
is trivial and there is a homotopy $\gamma_t$, $t\in[0,1]$, such that
$\gamma_0(x)=\Id$ for $x\in S_R^{p-1}$ and
$\gamma_1=\tilde\psi|S_R^{p-1}$. If $p=2$ and $\onm{deg}n=0$, then such
a homotopy $\gamma_t$ can be constructed by using the connection to
lift a homotopy between $n$ and the constant map $-\mbf k$. In any
case, choosing a smooth map $t(r)$ which is $0$ for $r\le R$ and $1$
for $r\ge R+1$ and setting $\psi(x)=\tilde\psi(x)(\gamma_{t(|x|)})^{-1}$
gives $\psi$ which lifts $n$ and is equal to the identity for $|x|\ge
R+1$. On the other hand, if $p=2$ and $\onm{deg}n\ne0$ then there is no
such $\psi$ since $\pi_2\bigl(SO(3)\bigr)=0$ and there would be a
homotopy of $\psi$ to a point which would project to a homotopy of $n$ to
a point.

For $p=2$ the degree nonzero sectors can only be obtained by reducing
at a nonconstant value for $\mu$, as opposed to $\mu=\mbf k$. 
The loss of translation invariance of $\mu$ leads to
serious complications that require more fully the reduction by stages theory. Towards the end of this Letter, I will attend
to this but for now assume $\mu=\mbf k$ and that $p\ne 2$, or $p=2$ and
$\onm{deg}n=0$.

The group $H=\mtl{SE}(p)=\bigl\{(U,a)\in SO(p)\times\mbb R^p\bigr\}$
also acts on $G$ by spatially rotating and translating its elements:
\begin{equation*}
(U,a)\cdot\psi=\psi\bigl((U,a)^{-1}x\bigr)
=\psi\bigl(U^{-1}(x-a)\bigr).
\end{equation*}
This action does not commute with the action of $G$, since
\begin{equation*}
(U,a)\cdot(A\cdot\psi)=(U,a)\cdot(A\psi)
=A\bigl(U^{-1}(x-a)\bigr)\psi\bigl(U^{-1}(x-a)\bigr)
\end{equation*}
whereas
\begin{equation*}
A\cdot\bigl((U,a)\cdot\psi\bigr)=A(x)\psi\bigl(U^{-1}(x-a)\bigr).
\end{equation*}
Since the actions of $G$ and $H$ do not commute they cannot be bound
together, using the direct product, into a single action. The
noncommutativity can be written, however, as
\begin{equation*}
(U,a)\cdot(A\cdot\psi)=\theta_{(U,a)}(A)\cdot\bigl((U,a)\cdot\psi\bigr)
\end{equation*}
as long as one defines $\theta:H\times G\rightarrow G$ by
\begin{equation*}
\theta\bigl((U,a),A\bigr)=A\bigl(U^{-1}(x-a)\bigr).
\end{equation*}
The assignment $(U,a)\mapsto\theta_{(U,a)}$ is a group morphism from
$H$ to $\onm{Aut} G$.

In general we have two groups $G$ and $H$ acting on $M$, the
actions do not commute, but there is a $\theta:H\rightarrow\onm{Aut}G$
such that for all $m\in M$, $g\in G$, and $h\in H$,
\begin{equation*}
h\cdot(g\cdot m)=\theta_h(g)\cdot(h\cdot m).
\end{equation*}
This situation is favorable to quotienting by stages since it is
a general situation where the action of $H$ on $M$ passes to
the quotient $M/G$. If one wants to define an action of $G\times H$ on $M$ by
$(g,h)\cdot m=g\cdot(h\cdot m)$, then
\begin{equation*}
(g_1,h_1)\cdot\bigl((g_2,h_2)\cdot m\bigr)=
g_1\cdot\bigl(h_1\cdot\bigl(g_2\cdot(h_2\cdot m)\bigr)\bigr)=
g_1\theta_{h_1}(g_2)\cdot\bigl((h_1h_2)\cdot m\bigr)
\end{equation*}
and the group product on $G\times H$ must be taken to be
the semidirect product
\begin{equation*}
(g_1,h_1)\cdot(g_2,h_2)=\bigl(g_1\theta_{h_1}(g_2),h_1h_2\bigr).
\end{equation*}
Thus, the two actions may be bound together but using the semidirect
product instead of the direct product. The semidirect product of $G$ and $H$
is commonly denoted $G\semi H$.

Returning to the Landau--Lifshitz equation, the momentum map for the
action of $H$ is easily calculated. Let $(\Omega,\dot a)\in\mfk h$,
where $\mfk h=\mfk{so}(p)\times\mbb R^p$ is the Lie algebra of $H$.
The infinitesimal generator of $(\Omega,\dot a)$ is, remembering to
use right translation,
\begin{equation*}\begin{align*}
&\dbyd\epsilon0\psi\bigl(\onm{exp}(-\Omega\epsilon)(x-\dot a\epsilon)
 \bigr)\psi^{-1}(x)\\
&\qquad\mbox{}=\dbyd\epsilon0\psi\bigl(\onm{exp}(-\Omega\epsilon)x
 \bigr)\psi^{-1}(x)
 +\dbyd\epsilon0\psi(x-\dot a\epsilon)\psi^{-1}(x)\\
&\qquad\mbox{}=\nabla^R_{-\Omega x-\dot a}\psi(x),
\end{align*}\end{equation*}
where $\nabla^R$, the right-hand gradient, is defined by
\begin{equation*}
(\nabla^R_b\psi)(x)=\dbyd\epsilon0\psi(x+\epsilon b)\psi(x)^{-1}.
\end{equation*}
I will use the inner product
$\Omega_1\cdot\Omega_2=-\frac12\onm{trace}(\Omega_1\Omega_2)$ for
$\Omega_1,\Omega_2\in\mfk{so}(p)$ and the notation $v\wedge
w=v\otimes w-w\otimes v=vw^t-wv^t$ for $v,w\in\mbb R^p$, 
so that
\begin{equation*}
x\cdot\Omega y=\langle x\wedge y,\Omega\rangle.
\end{equation*}
As well, I
will set
\begin{equation*}
(\mu\cdot\nabla^R)\psi=(\mu\cdot\nabla^R_{\mbf{e_1}}\psi,\cdots,
\mu\cdot\nabla^R_{\mbf{e_p}}\psi),
\end{equation*}
where $\mbf{e_1},\cdots,\mbf{e_p}$ is the standard basis of $\mbb R^p$, so that
\begin{equation*}
v\cdot\nabla^R_w\psi=w\cdot(v\cdot\nabla^R\psi),\quad v,w\in\mbb R^p.
\end{equation*}
Then
\begin{equation*}\begin{align*}
J^H_{\Omega,\dot a}(\psi,\mu)&=\int\mu\cdot\nabla^R_{-\Omega x
 -\dot a}\psi(x)\\
&=(\Omega,\dot a)\cdot\left(\int x\wedge(\mu\cdot\nabla^R)\psi,
 -\int(\mu\cdot\nabla^R)\psi\right),
\end{align*}\end{equation*}
so the $H$ momentum map is
\begin{equation*}
J^H(\psi,\mu)=\left(\int x\wedge(\mu\cdot\nabla^R)\psi,-\int(\mu\cdot
\nabla^R)\psi\right).
\end{equation*}
The reduction of the symmetry represented by the $H$ action
is possible after the following two observations:
\begin{itemize}
\item\emph{$(J^G)^{-1}(\mu)$ is an invariant subset under the
action of $H$.} Indeed,
\begin{equation*}
J^G\bigl((U,a)\cdot(\psi,-\psi\mbf k)\bigr)=
-\psi\bigl((U,a)^{-1}x\bigr)^{-1}\psi\bigl((U,a)^{-1}x\bigr)\mbf k=\mbf k.
\end{equation*}
\item\emph{$J^H$ is $G_\mu$ invariant for $p\ge 2$.} 
Assume $A\in G_\mu$. Then
using the (easily verified) identity
\begin{equation*}
(\mu\cdot\nabla^R)(\psi\phi)=(\mu\cdot\nabla^R)\psi+\bigl((\psi^{-1}
\mu)\cdot\nabla^R\bigr)\phi,
\end{equation*}
one has
\begin{equation*}\begin{split}
&J^H\bigl(A\cdot(\psi,-\psi\mbf k)\bigr)\\
&\qquad\mbox{}=\biggl(\int x\wedge\bigl((\psi\mbf k\bigr)\cdot\nabla^R)
 (\psi A^{-1}),
 -\int\bigl((\psi\mbf k)\cdot\nabla^R\bigr)(\psi A^{-1})\biggr)\\
&\qquad\mbox{}=J^H(\psi,-\psi\mbf k)+
\biggl(\int x\wedge\bigl((\mbf k\cdot\nabla^R)A^{-1}\bigr),
 -\int(\mbf k\cdot\nabla^R) A^{-1}\biggr).\\
\end{split}\end{equation*}
Since $A(x)\mbf k=\mbf k$, $A(x)=\onm{exp}\bigl(\alpha(x)\mbf k\bigr)$
for some smooth real valued function $\alpha$ which is an integral multiple of $2\pi$ at infinity, one has
\begin{equation*}
(\mbf k\cdot\nabla^R)A^{-1}=-\nabla\alpha,
\end{equation*}
so
\begin{equation*}
J^H\bigl(A\cdot(\psi,-\psi\mbf k)\bigr)=J^H(\psi,-\psi\mbf k)+
\biggl(-\int x\wedge\nabla\alpha,
\int\nabla\alpha\biggr).
\end{equation*}
In the case $p\ge2$ the space at $\infty$ is connected so $\alpha$ is
constant there and hence both integrals vanish.
\end{itemize}

The nontrivial part of the above is the second item about the
invariance of $J^H$; the first item is a specialty of the particular
momentum $\mu$. Yet this kind of thing occurs generally. I
will use the notation $\onm{CoAd}_g\mu=\onm{Ad}^*_{g^{-1}}\mu$.

\begin{lem}\label{lm1}Let $J^{\tilde G}:P\rightarrow\tilde{\mfk g}^*$
be an equivariant momentum map and let $G$ be a normal subgroup of $\tilde G$. Define $J^G=i_{\mfk
g}^* J^{\tilde G}$, so that $J^G$ is an equivariant momentum map
for the action of $G$, let $\mu\in\mfk g^*$, define the
subgroup $H_\mu$ of $\tilde G$ by
\begin{equation*}
H_\mu=\bset{g\in\tilde G}{\mbox{$J^G(gp)=\mu$ for all $p\in (J^G)^{-1}(\mu)$}},
\end{equation*}
and define $J^{H_\mu}=i_{\mfk h_\mu}^*J^{\tilde G}$.
Suppose that $G_\mu$ is connected.
Then $G_\mu=G\cap H_\mu$, $G_\mu$ is a normal subgroup of $H_\mu$, and $J^{H_\mu}$ is $G_\mu$ invariant on $(J^G)^{-1}(\mu)$. 
\end{lem}
\begin{proof}
For $g\in G$ the condition $J^G(gp)=\mu$ for all $p\in (J^G)^{-1}(\mu)$
is exactly $\onm{CoAd}_g\mu=\mu$ so $G_\mu=G\cap H_\mu$. Then $G_\mu$
normal in $H_\mu$ follows directly from $G$ normal in $\tilde G$.

It is required to prove that if $p\in P$ is such that $J^G(p)=\mu$, if
$g\in G_\mu$, and if $\xi\in\mfk h_\mu$, then $J_\xi^{\tilde
G}(gp)=J_\xi^{\tilde G}(p)$. Setting $\tilde\mu=J^{\tilde G}(p)$ and, since
$G_\mu$ is connected, this is equivalent to
\begin{equation*}
\dbyd t0J_\xi^{\tilde G}\bigl(\onm{exp}(\eta t)p\bigr)=
 \langle J^{\tilde G}(p),[\xi,\eta]\rangle= \langle\tilde\mu,
 [\xi,\eta]\rangle=0
\end{equation*}
for all $\eta\in\mfk g_\mu$. But if $\xi\in\mfk h_\mu$ and
$\eta\in\mfk g_\mu$ then by the definition of $H_\mu$
\begin{equation*}
\langle\onm{CoAd}_{\onm{exp}(\xi t)}\tilde\mu,\eta\rangle=\langle\mu,
 \eta\rangle
\end{equation*}
so differentiation in $t$ gives
$\langle\tilde\mu,[\xi,\eta]\rangle=0$. Then $[\xi,\eta]\in\mfk
g_\mu\subseteq\mfk g$ since $\eta\in\mfk g_\mu$ and $\xi\in\mfk
h_{\tilde\mu}$ and $G_\mu$ is normal in $H_\mu$, so
$\langle\tilde\mu,[\xi,\eta]\rangle=\langle\mu,[\xi,\eta]\rangle=
\langle\onm{coad}_\eta\mu,\xi\rangle=0$.\eop
\end{proof}

As is easily verified, since the subgroup $G$ is assumed to be normal,
the coadjoint action of $\tilde G$ restricts to $\mfk g$ and the
subgroup $H_\mu$ is the isotropy group of $\mu$ with respect to this
action.

The point of Lemma~\ref{lm1} is that when reducing first by the subgroup
$G$ at the momentum value $\mu$ to obtain $(J^G)^{-1}(\mu)/G_\mu$, the
residual symmetry after the reduction would intuitively be those
elements $\tilde g\in \tilde G$ which map $(J^G)^{-1}(\mu)$ to itself;
this is exactly the group $H_\mu$. Then, in order to pass the
`residual momentum' to the quotient it is required that $J^{H_\mu}$
be $G_\mu$ invariant, and the Lemma gives sufficient conditions for
this. As for the relevance to the Landau--Lifshitz equation, the
Euclidean invariance of momentum $\mu=\mbf k$ gives $H_\mu=G_\mu\semi
H$ (that is rotation fields fixing $\mbf k$ paired with all Euclidean
motions), $G$ is normal in the $G\semi H$, and so the Lemma explains
the invariance of $J^H$ under $G_\mu$ for $m\ge2$. Also, for
dimension $q=1$, the group $G_\mu$ is not connected, having components
in one-to-one correspondence with $\mbb Z$, and so the invariance of
$J^H$ fails when $q=1$, just when a hypothesis of the Lemma fails.

Obtaining a formula for the momentum for the Landau--Lifshitz equation
means actually passing the momentum map
$J^H$ to the quotient space $\{n\}$.
One way to do this is, given some $n$, to find $\psi_n$ such that
$\psi_n\mbf k=-n$ and then calculate $J^H(\psi_n)$, and one such
$\psi_n$ is
\begin{equation*}
\psi_n=\onm{exp}\left(\frac{-\onm{arccos}(-\mbf k\cdot n)}{|\mbf k\times n|}
\mbf k\times n\right).
\end{equation*}
A long calculation gives
\begin{equation*}
(n\cdot\nabla^R)\psi_n=\frac{-1}{1-\mbf k\cdot n}
\bigl((\mbf k\times n)\cdot\nabla\bigr)n,
\end{equation*}
so the momentum map on the reduced space is
\begin{equation*}
\bar J^H(n)=\left(\int x\wedge\frac{1}{1-\mbf k\cdot n}\bigl((\mbf k\times n)
\cdot\nabla\bigr)n,\int\frac{-1}{1-\mbf k\cdot n}\bigl((\mbf k\times n)\cdot
\nabla\bigr)n\right).
\end{equation*}
This is somewhat unsatisfactory due to the singularity at $n=\mbf k$,
which is due to the singularity of $\psi_n$ at $n=\mbf k$. A concern
might be that the singularities in the choice of $\psi_n$ might give a
value for the momentum that differs from nonsingular choices but the
singularities are tame enough not to contribute to the integrals
making up the momenta.

There remains the problem of the degree nonzero sectors in the case
$p=2$. Hitting the degree $m$ sector is not a problem: it is the
reduced space $(J^G)^{-1}(\mu)/G_\mu$ for some $\mu\in\mfk g$ such that
$\onm{deg}\mu=m$ and $|\mu(x)|=1$ for all $x\in\mbb R^2$. The
difficulty is that the subgroup $H$, which ought to
represent Euclidean symmetries, does not act on
$(J^G)^{-1}(\mu)$, since
\begin{equation*}
\begin{split}
J^G\bigl((U,a)(\psi,-\psi\mu)\bigr)
 &=J^G\Bigl(\psi\bigl(U^{-1}(x-a)\bigr),-\psi\bigl(U^{-1}(x-a)\bigr)
 \mu\bigl(U^{-1}(x-a)\bigr)\Bigr)\\
 &=\mu\bigl(U^{-1}(x-a)\bigr)\ne\mu(x).
\end{split}
\end{equation*}
This is a fundamental obstruction and it is best to have recourse to
the general theory. An appropriate general context is that
of~Lemma~\ref{lm1}. Since $G_\mu$ is normal in $H_\mu$, the quotient
group $\bar H_\mu=H_\mu/G_\mu$ acts on $P_\mu=(J^G)^{-1}(\mu)/G_\mu$ and the
obvious adjustment is to use the group $\bar H_\mu$ to represent the
`residual' symmetries on~$P_\mu$. Since $J^{H_\mu}$ is $G_\mu$ invariant
it passes to the quotient, but it is not clearly a momentum map for
the action of the group $\bar H_\mu$ since it has values in $\mfk
h_\mu^*$ rather than $\bar\mfk h_\mu^*$.
The values of $J^{H_\mu}$ are in $\bar\mfk h_\mu^*=(\mfk h_\mu/\mfk g_\mu)^*$ if and only if they annihilate
$\mfk g_\mu$. This cannot be expected since if $p\in
(J^G)^{-1}(\mu)$ and $\xi\in\mfk g_\mu$ then
\begin{equation*}
\langle J^{H_\mu}(p),\xi\rangle=\langle\mu,\xi\rangle,
\end{equation*}
which vanishes for all $\xi\in\mfk g_\mu$ if and only if $\mu$
annihilates $\mfk g_\mu$, an unusual circumstance. For example, if
$\mu\ne0$, and in the presence of an $\onm{Ad}$ invariant metric with
the usual identification of $\mfk g_\mu^*$ and $\mfk g_\mu$, one would
have $\mu\in\mfk g_\mu$ and $\mu$ would not annihilate itself. Also if
$\mfk g$ is Abelian then $\mu$ annihilates $\mfk g_\mu$ if and only if
$\mu$ is zero. The obvious adjustment here is to subtract from
$J^{H_\mu}$ any extension, say $\hat\mu$, of $\mu|\mfk g_\mu$ to $\mfk
h_\mu$. Then $J^{H_\mu}-\hat\mu$ has values in $\bar\mfk h_\mu$ and
is $G_\mu$ invariant on $(J^{G})^{-1}(\mu)$, and so defines a map
$J^{\bar H_\mu}:P_\mu\rightarrow\bar\mfk h_\mu^*$. This is a momentum
map for the action of $\bar H_\mu$ of $P_\mu$. Indeed, given
$\xi\in\bar\mfk h_\mu=\mfk h_\mu/\mfk g_\mu$ there is a $\hat\xi\in\mfk h_\mu$ projecting to
$\xi$, and the flow of $J^{H_\mu}_{\hat\xi}-\hat\mu$, which is
multiplication by $\onm{exp}(\hat\xi t)$, projects to the flow of
$J^{\bar H_\mu}_\xi$. By construction the projection of
multiplication by $\onm{exp}(\hat\xi t)$ is multiplication by
$\onm{exp}(\xi t)$, so $J^{\bar H_\mu}$ generates the action of $\bar
H_\mu$, which is what is required for $J^{\bar H_\mu}$ to be a
momentum map. This reduction procedure is identical to that found in~\cite{MarsdenJEMisiolekGPerlmutterMRatiuTS-1999}.

What is interesting about all this is that the momentum map
$J^{\bar H_\mu}$ \emph{need not be equivariant, even though it has
been reduced from an equivariant momentum map.} 
This effect is also observed in~\cite{MarsdenJEMisiolekGPerlmutterM-1998.1}
and~\cite{MarsdenJEMisiolekGPerlmutterMRatiuTS-1999}.
The lack of
equivariance is characterized (see~\cite{MarsdenJERatiuTS-1994.1}) by
the cocycle $\sigma^{\bar H_\mu}=J^{\bar H_\mu}(hp)-\onm{CoAd}_hJ^{\bar
H_\mu}(p)$, $h\in\bar H_\mu$. The value of $p$ is irrelevant as long as $P_\mu$
is connected. To calculate this cocycle,
given $p\in P_\mu$, choose $\hat p\in (J^{G})^{-1}(\mu)$ projecting to
$p$, and given $h\in\bar H_\mu$ choose $\hat h\in H_\mu$ projecting to
$h$. Then
\begin{equation*}
\sigma^{\bar H_\mu}(h)
=J^{H_\mu}(\hat h\hat p)
 -\hat\mu-\onm{CoAd}_{\hat h}\bigl(J^{H_\mu}(\hat p)-\hat\mu\bigr)
=\onm{CoAd}_{\hat h}\hat\mu-\hat\mu.
\end{equation*}
The values of $\sigma^{\bar H_\mu}$ annihilate $\mfk g_\mu$ and are to
be regarded as lying in $\bar\mfk h_\mu^*$. The infinitesimal version
of this cocycle is the antisymmetric, bilinear cocycle $\Sigma^{\bar H_\mu}$
on $\bar\mfk h_\mu$
defined by
\begin{equation*}
\Sigma^{\bar H_\mu}(\xi,\eta)
 =\dbyd\epsilon0\bigl\langle\sigma^{\bar H_\mu}
 \bigl(\onm{exp}(\eta\epsilon)\bigr),
 \xi\bigr\rangle.
\end{equation*}
Given $\xi,\eta\in\bar\mfk h_\mu$, and choosing $\hat\xi,\hat\eta\in\mfk h_\mu$
projecting to $\xi,\eta$ respectively, it is immediate that
\begin{equation*}
\Sigma^{\bar H_\mu}(\xi,\eta)=\langle\hat\mu,[\hat\xi,\hat\eta]\rangle.
\end{equation*}
Nonequivariance of the momentum map is reflected in the commutation relations
\begin{equation*}
\{J^{\bar H_\mu}_\xi,J^{\bar H_\mu}_\eta\}=J^{\bar H_\mu}_{[\xi,\eta]}
 -\Sigma^{\bar H_\mu}(\xi,\eta)
\end{equation*}
which differ from the usual ones where no cocycle appears.

The reduction just described can be worked out in
detail for the Landau--Lifshitz equation. Fix some $\mu\in\mfk g$ which
is not necessarily constant but $|\mu(x)|=1$; when $p=2$ this $\mu$
could have nonzero degree. The subgroup $H_\mu$ is found by
finding the subgroup under which $(J^G)^{-1}(\mu)$ is mapped into itself.
Since
\begin{equation*}
\begin{split}
&J^G\Bigl(\bigl(A,(U,a)\bigr)(\psi,-\psi\mu)\Bigr)\\
&\qquad\mbox{}=J^G\Bigl(\psi\bigl(U(x)^{-1}(x-a)\bigr)A(x)^{-1},
 -\psi\bigl(U(x)^{-1}(x-a)\bigr)\mu\bigl(U(x)^{-1}(x-a)\bigr)\Bigr)\\
&\qquad\mbox{}=-A(x)\psi\bigl(U(x)^{-1}(x-a)\bigr)^{-1}
 \times-\psi\bigl(U(x)^{-1}(x-a)\bigr)\mu\bigl(U(x)^{-1}(x-a)\bigr)\\
&\qquad\mbox{}=A(x)\mu\bigl(U(x)^{-1}(x-a)\bigr)
\end{split}
\end{equation*}
it follows that
\begin{equation*}
H_\mu=\bset{\bigl(A,(U,a)\bigr)}{A(x)\mu\bigl(U^{-1}(x-a)\bigr)=\mu(x)}.
\end{equation*}
The projection onto the second factor is a quotient map for the left
action of $G_\mu$ on $H_\mu$, so $\bar H_\mu= H_\mu/G_\mu\cong H$ and the
action of $\bar H_\mu$ on the space of Landau--Lifshitz fields 
$P_\mu=\{n\}$ is
\begin{equation*}
(U,a)n=n\bigl(U(x)^{-1}(x-a)\bigr),
\end{equation*}
as expected. 

To calculate the infinitesimal cocycle first extend $\mu$ to
$\mfk h_\mu$. The extension that comes to mind is
\begin{equation*}
\bigl\langle\hat\mu,\bigl(\xi,(\Omega,\dot a)\bigr)\bigr\rangle
 =\int\mu(x)\cdot\xi(x),
\end{equation*}
or just $\iota_1^*\mu$ where $\iota_1:G\rightarrow G\semi H$ is the
inclusion. The subgroup $H_\mu$ is defined by the
condition $A(x)\mu\bigl(U^{-1}(x-a)\bigr)=\mu(x)$, so its Lie algebra $\mfk
h_\mu$ is obtained by differentiating this condition at the identity.
Given $\bigl(\xi,(\Omega,\dot a)\bigr)\in\mfk{se}(p)$, this gives
\begin{equation*}
\dbyd\epsilon0\onm{exp}\bigl(\xi(x)\bigr)\mu\bigl(
 \onm{exp}(\Omega\epsilon)(x-\dot a\epsilon)\bigr)
=\xi(x)\times\mu(x)+\nabla_{-\Omega x-\dot a}\mu(x)=0,
\end{equation*}
so given an $(\Omega,\dot a)\in\bar H_\mu\cong\mfk{se}(p)$ one finds a
$\bigl(\xi,(\Omega,\dot a)\bigr)\in\mfk{se}(p)\in\mfk h_\mu$ projecting to this
by solving
\begin{equation*}
\xi(x)\times\mu(x)=\nabla_{\Omega x+\dot a}\mu(x).
\end{equation*}
One solution is
\begin{equation*}
\xi(x)=\mu(x)\times\nabla_{\Omega x+\dot a}\mu(x)
\end{equation*}
so that one can set
\begin{equation*}
(\Omega,\dot a)^\wedge=\bigl(\mu(x)\times\nabla_{\Omega x+\dot a}\mu(x),
(\Omega,\dot a)\bigr).
\end{equation*}
The Lie bracket of $G\semi H$ is
\begin{equation*}
\begin{split}
&\bigl[\bigl(\xi_1,(\Omega_1,\dot a_1)\bigr),
 \bigl(\xi_2,(\Omega_2,\dot a_2)\bigr)\bigr]\\
&\qquad\mbox{}=\bigl(\xi_1(x)\times\xi_2(x)-\nabla_{\Omega_1x+\dot a_1}\xi_2(x)+
 \nabla_{\Omega_2x+\dot a_2}\xi_1(x),[(\Omega_1,a_1),(\Omega_2,a_2)]\bigr)
\end{split}
\end{equation*}
so the cocycle is (many of the terms vanish because they are obviously
perpendicular to $\mu$ and so vanish when paired with $\mu$ in the
outermost dot product)
\begin{equation*}
\Sigma^{\bar H_\mu}\bigl((\Omega_1,\dot a_1),(\Omega_2,\dot a_2)\bigr)=
 -\int\mu(x)\cdot\nabla_{\Omega_1x+\dot a_1}\mu(x)\times
 \nabla_{\Omega_2x+\dot a_2}\mu(x).
\end{equation*}
When $\mu=\mbf k$ this cocycle is zero. If $\mu$ is spherically
symmetric, then $\Sigma^{\bar H_\mu}$ simplifies to
\begin{equation*}
\Sigma^{\bar H_\mu}\bigl((\Omega_1,\dot a_1),(\Omega_2,\dot a_2)\bigr)=
-\int\mu(x)\cdot\nabla_{\dot a_1}\mu(x)\times
 \nabla_{\dot a_2}\mu(x),
\end{equation*}
which in the the case $p=2$ is
\begin{equation*}
\Sigma^{\bar H_\mu}\bigl((\Omega_1,\dot a_1),(\Omega_2,\dot a_2)\bigr)=
-4\pi\omega_0(\dot a_1,\dot a_2)\onm{deg}\mu,
\end{equation*}
where 
\begin{equation*}\omega_0(a_1,a_2)=a_1^tJa_2, \quad
J=\left[\begin{array}{cc}0&1\\-1&0\end{array}\right].
\end{equation*}
A similar $\mfk{se}(2)$ cocycle appears in the planar point vortex
system~\cite{PatrickGW-1999.1}.

In order for the further reduction of the Landau--Lifshitz system by
the symmetry $\bar H_\mu$ to be symplectomorphic to the original
system on $T^*G$ reduced by $G\semi H$, an additional `stages
hypothesis'~(\cite{MarsdenJEMisiolekGPerlmutterMRatiuTS-1999}) must
be satisfied. It is sufficient that $\mfk g+\mfk h_\mu=\mfk
g\semi\mfk h$, a fact which is obvious from the calculation of $\mfk
h_\mu$ above.

The momentum map $J_{\Omega,\dot a}^{\bar H_\mu}(n)$ is to be
calculated by finding $J_{(\Omega,\dot
a)^\wedge}^{H_\mu}(\psi,-\psi\mu)$ where $\psi$ such that $\psi\mu=n$.
With the chosen extension $\hat\mu$ of $\mu$ the difference
$J^{H_\mu}-\hat\mu$ annihilates all of the first factor of $\mfk
g\semi\mfk h$, so this is equivalent to calculating $J_{\Omega,\dot
a}^H(\psi,-\psi\mu)$. The result is that the reduced momentum is
calculated in the same way as the case $\mu=\mbf k$, which is just to
find $\psi$ such that $\psi\mu=n$ and then calculate
$J^H(\psi,-\psi\mu)$. If one then takes $n=\mu$, the momentum
is zero since one can take $\psi(x)=\Id$. Thus, the reduced
momentum calculated this way is the unique momentum that vanishes at $\mu$.

As has already been mentioned, from~\cite{CooperNR-1999.1}
and~\cite{PapanicolaouNTomarasTN-1991.1} the conserved quantities
\begin{equation*}
P=2\pi\int x\times\Omega,\qquad\Omega_\alpha=\frac1{8\pi}
 \epsilon_{\alpha\beta\gamma}n\cdot(\nabla_\beta n\times
 \nabla_\gamma n),
\end{equation*}
generate translations in the case $p=3$. 
Actually, working out the $x$-component of this where, $\mbb
R^3=\{(x,y,z)\}$, easily gives 
\begin{equation*}
P_x=\frac12\int n\cdot\left(\frac{\partial n}{\partial x}\times\left(
y\frac{\partial n}{\partial y}+z\frac{\partial n}{\partial z}\right)\right),
\end{equation*}
and leads directly to the generalization (summation over repeated indices)
\begin{equation*}
P_i=\frac{1}{p-1}\int n\cdot\left(\frac{\partial n}{\partial x^i}\times
x^j\frac{\partial n}{\partial x^j}\right).
\end{equation*}
It is easily verified that this generates translations for any value of $p\ne 1$.
If a vector density $\mcl P$ is defined by
\begin{equation*}
\mcl P_i=\frac1{p-1}n\cdot\left(\frac{\partial n}{\partial x^i}\times
x^j\frac{\partial n}{\partial x^j}\right),
\end{equation*}
then, in view of the expression for $J^H$, one might guess that
$-\Omega\cdot(x\wedge\mcl P)=\Omega x\cdot\mcl P$ is a density
generating the rotations corresponding to $\Omega$. This is also
easily verified, and so leads to the momentum map
\begin{equation*}
\tilde J^{\bar H_\mu}=
\left(-\int x\wedge\mcl P,\int\mcl P\right).
\end{equation*}
When $n$ is spherically symmetric, the density $\mcl P$ is zero since
\begin{equation*}\begin{split}
\mcl (p-1)P_i&=n\cdot\left(\frac{\partial n}{\partial x^i}\times 
 x^j\frac{\partial n}
 {\partial x^j}\right)\\
&=n\cdot\left(\frac{\partial n}{\partial r}\frac{\partial r}
 {\partial x^i}
 \times x^j\frac{\partial n}{\partial r}\frac{\partial r}{\partial x^j}\right)
= n\cdot\left(\frac{\partial n}{\partial r}\frac{x^i}{r}\times
 r\frac{\partial n}{\partial r}\right) = 0.
\end{split}
\end{equation*}
Thus, $J^{\bar H_\mu}=\tilde J^{\bar H_\mu}$, and the momentum
of~Papanicolaou and~Tomaras~\cite{PapanicolaouNTomarasTN-1991.1} is
seen to be the result of semidirect product reduction. In particular,
their observation that in the case $p=2$ one has
$\{P_x,P_y\}(n)=4\pi\onm{deg}n$ follows from the commutation relations
and the cocycle, as follows:
\begin{equation*}
\begin{split}
\{P_x,P_y\}(n)&=\{J_{\mbf i}^{\bar H_\mu},J_{\mbf j}^{\bar H_\mu}\}(n)\\
&=J_{[\mbf i,\mbf j]}^{\bar H_\mu}(n)-\Sigma^{\bar H_\mu}(\mbf i,\mbf j)
=4\pi\omega_0(\mbf i,\mbf j)\deg\mu
=4\pi\deg n.
\end{split}
\end{equation*}

\frenchspacing
\bibliographystyle{plain}

\end{document}